\documentclass[twoside]{article}
\usepackage{landau}
\usepackage[dvips]{epsfig}
\theoremstyle{definition}

\theoremstyle{remark}

\numberwithin{equation}{section}

\newcommand{\eps}{\varepsilon}

\newcommand{\del}{\partial}
\newcommand{\R}{\mathbb{R}}
\newcommand{\e}{{\varepsilon}}
\newcommand{\bla}{{\boldsymbol \lambda}}
\newcommand{\x}{\mathbf{x}}
\newcommand{\vv}{\mathbf{v}}
\renewcommand{\r}{\rho}
\renewcommand{\xi}{{s}}
\renewcommand{\d}{\partial}
\newcommand{\G}{\Gamma}
\renewcommand{\k}{\kappa}
\newcommand{\rr}{{\bf r}}

\def\shorttitle{Weak Continuity and Compactness for Nonlinear PDEs}
\def\shortauthor{Gui-Qiang G. Chen}

\begin{document}

\title
{\Large \bf \boldmath Weak Continuity and Compactness for Nonlinear Partial Differential
Equations\footnote{The research of Gui-Qiang Chen was supported in part by
the UK EPSRC Science and Innovation
Award to the Oxford Centre for Nonlinear PDE (EP/E035027/1),
the UK EPSRC Award to the EPSRC Centre for Doctoral Training
in PDEs (EP/L015811/1),
the NSFC under a joint project Grant 10728101, and
the Royal Society-Wolfson Research Merit Award (UK).}}

\author{\large Gui-Qiang G. Chen
    \\ \normalsize\emph{School of Mathematical Sciences, Fudan University,}
    \\ \normalsize\emph{Shanghai 200433, China.}
    \\ \normalsize\emph{Academy of Mathematics and Systems Science,}
    \\ \normalsize\emph{Chinese Academy of Sciences, Beijing 100190, China.}
        \\ \normalsize\emph{Mathematical Institute, University of Oxford,}
    \\ \normalsize\emph{Oxford, OX2 6GG, UK.}
        \\ \normalsize\emph{E-mail: Gui-Qiang.Chen@maths.ox.ac.uk}
}

\date{\vspace{-12mm}}

\maketitle

\thispagestyle{first}

\bigskip
% ----------------------------------------------------------------
\begin{abstract}
This paper presents several examples of fundamental
problems involving weak continuity and compactness for
nonlinear partial differential equations,
in which compensated compactness and related ideas
have played a significant role.
The compactness and convergence of vanishing viscosity
solutions for nonlinear hyperbolic conservation laws
are first analyzed,
including the inviscid limit from the
Navier-Stokes equations to the Euler equations for homentropic flow,
the vanishing viscosity method to construct the global
spherically symmetric solutions to the multidimensional
compressible Euler equations,
and the sonic-subsonic limit of solutions of the full Euler
equations for multidimensional steady compressible fluids.
The weak continuity and rigidity
of the Gauss-Codazzi-Ricci system and corresponding
isometric embeddings in
differential geometry are revealed.
Further references are also provided for some recent
developments on the weak continuity and compactness
for nonlinear partial differential equations.
\end{abstract}

% ----------------------------------------------------------------
\section{Introduction}

Nonlinear partial differential equations (PDEs) can be written as
the following general form:
\begin{equation}\label{1.1}
\mathcal{N}[U] = 0,
\end{equation}
where $\mathcal{N}[\cdot]$ is a nonlinear mapping, and $U$ is an unknown (scalar or vector) function that
is called a solution if $U$ solves \eqref{1.1}.

\smallskip
Two of the fundamental issues for nonlinear PDEs \eqref{1.1} are the following:
\begin{enumerate}
\renewcommand{\theenumi}{\roman{enumi}}
\item[(i)] {\it Weak Continuity and Rigidity}: $\,$ Let $\{U^\eps\}_{\eps>0}$
be a sequence of exact solutions satisfying
\begin{equation}\label{1.2}
\begin{cases}
\mathcal{N}[U^\varepsilon]=0,\\[0.5mm]
U^\varepsilon \rightharpoonup U  \qquad \text{in some topology as $\e \to 0$}.
\end{cases}
\end{equation}

{\bf Issue 1}: {\it Does the limit function $U$ satisfy
\begin{equation}\label{1.3}
 \mathcal{N}[U]=0,
\end{equation}
or
\begin{equation}\label{1.4}
\tilde{\mathcal{N}}[U]=0
\end{equation}
for  a different nonlinear mapping $\tilde{\mathcal{N}}[\cdot]$
associated with
the original nonlinear mapping $\mathcal{N}[\cdot]$ and the solution sequence
$\{U^\eps\}_{\eps>0}$?}

\smallskip
Such an issue arises in rigidity problems in geometry, mechanics,
among others.

\item[(ii)] {\it Compactness and Convergence}: $\,$ Let $\{U^\varepsilon\}_{\varepsilon>0}$
be a sequence of approximate or multiscale solutions satisfying
\begin{equation}\label{1.5}
\begin{cases}
\mathcal{N}^\varepsilon[U^\varepsilon]=0,\\[0.5mm]
U^\varepsilon \rightharpoonup U  \qquad \text{in some topology as $\e\to 0$}.
\end{cases}
\end{equation}

{\bf Issue 2}: {\it Does the limit function $U$ satisfy
\eqref{1.3}, or \eqref{1.4} for a different nonlinear mapping
$\tilde{\mathcal{N}}[\cdot]$
associated with
the nonlinear mappings $\mathcal{N}^\eps[\cdot]$ and the solution sequence
$\{U^\eps\}_{\eps>0}$?}

\smallskip
This issue arises in the viscosity methods,
relaxation methods, numerical methods, as well as
problems for homogenization, hydrodynamic limits,
search for effective equations, among others.
\end{enumerate}

This paper presents several examples of these fundamental problems
involving weak continuity and compactness for nonlinear PDEs,
in which compensated compactness and related ideas,
developed by Luc Tartar \cite{Tartar}--\cite{Tartar87}
and Fran\c{c}ois Murat \cite{Murat}--\cite{Murat87},
have played a significant role; also see Tartar \cite{Tartar09}.
In particular, in Section 2, we first analyze the compactness
and convergence of vanishing viscosity solutions
to hyperbolic conservation laws.
In Section 3, we reveal the weak continuity and rigidity
of the Gauss-Codazzi-Ricci system and corresponding isometric embeddings
in differential geometry.
Further references are also provided for some recent developments
on the weak continuity and compactness for nonlinear PDEs.
We finally remark that, as we will see in Sections 3--4,
many fundamental problems in this direction
are still open, which require further new mathematical ideas, techniques,
and approaches that deserve our special attention.

\section{Compactness and Convergence of Vanishing Viscosity Solutions
to Hyperbolic Conservation Laws}

\smallskip
Consider the following one-dimensional nonlinear hyperbolic conservation
laws with form:
\begin{equation}\label{2.1}
\del_t U + \del_x F(U) =0, \qquad U \in \R^N,
\end{equation}
where $F: \R^N\to \R^N$ is a nonlinear mapping so that all the eigenvalues of $\nabla_UF(U)$ are real.

To solve these nonlinear PDEs, one of the important approaches
is the viscosity method for which one
honors the {\it physical} or designs an {\it artificial} $N\times N$ matrix function:
\begin{equation}\label{2.2}
D: \R^N\to M^{N\times N}, \qquad D(U)\ge 0,
\end{equation}
so that
\begin{enumerate}
\renewcommand{\theenumi}{\roman{enumi}}
\item[(i)] $\del_t U + \del_x F(U)
=\varepsilon \partial_x\big(D(U)\del_xU\big)$
admits a global solution {$U^\e(t,x)$} for each fixed {$\e>0$};

\item[(ii)] $U^\varepsilon(t,x)\to U(t,x)$  in some topology as $\e\to 0$,
and $U(t,x)$ is an entropy solution.
\end{enumerate}
This method for the multidimensional case can be analogously formulated.

\medskip
The idea of the vanishing viscosity method
originates the philosophy of regarding the inviscid gas
as the limit of viscous gases, which can date back in the 19th century,
including
the work by Stokes (1848),
Rankine (1870), Hugoniot (1889),
Rayleigh (1910), Taylor (1910), Weyl (1949), among others; also
see Dafermos \cite{Dafermos} and the references cited therein.
This idea has played an essential role in developing the mathematical theory of
hyperbolic conservation laws (such as discontinuous solutions, entropy conditions,
existence, uniqueness, and solution behavior),
as well as numerical methods and related applications
(such as shock capturing, upwind, and kinetic schemes).
This method becomes increasingly important, especially
for understanding the recently observed non-uniqueness
phenomena for the weak solutions satisfying the entropy equality
for the multidimensional Euler equations
({\it cf.} \cite{DeLellis-Szekelyhidi-1,DeLellis-Szekelyhidi-2}).
On the other hand, the realization of this method is truly challenging
in mathematics, since it involves several fundamental difficulties in
analysis, including singular limits, nonlinearity, discontinuity, singularity, oscillation,
cavitation, and concentration.

%=======================================================================

\subsection{Compactness and Convergence via $BV$--Estimates}

This compactness framework is based on the compactness theorem in BV,
which is a sufficient framework to ensure the strong compactness and convergence
of exact/approximate solutions.
On the other hand, achieving the BV--estimates of exact/approximate solutions is
usually very challenging for the nonlinear systems, even though it is relatively
easier for the scalar case.

\subsubsection{Scalar conservation laws}

Consider the Cauchy problem for scalar conservation laws ($N=1$):
\begin{equation}\label{2.3}
\partial_t U+\partial_x F(U)=\varepsilon \partial_{xx} U,
\end{equation}
with the initial data $U|_{t=0}=U_0 \in BV\cap L^\infty(\R)$.
It can be shown that there exists $C$ independent of $\e$ such that
the viscous solutions $U^\eps=U^\eps(t,x)$  of \eqref{2.3} satisfy
\begin{enumerate}
\renewcommand{\theenumi}{\roman{enumi}}
\item[(i)] {\it Maximum principle}: $\,\,\,\|U^\e\|_{L^\infty} \le C$;

\item[(ii)] {\it BV--estimate}: $\,\,\,\|\partial_x U^\e\|_{L^1}+\|\partial_t U^\e\|_{L^1}\le C$.
\end{enumerate}
See Hopf \cite{Hopf}, Oleinik \cite{Oleinik}, and Lax \cite{Lax73} for the one-dimensional case,
and Vol'pert \cite{Volpert}
and Kruzhkov \cite{Kruzkov} for the multidimensional case.

\smallskip
One of the approaches to achieve the BV--estimate
is due to Vol'pert \cite{Volpert}, which yields
\begin{eqnarray*}
&&\partial_t (|\partial_x U^\e|) +\partial_x (F'(U^\e)|\partial_x U^\e|)
  \le \varepsilon \partial_{xx} (|\partial_x U^\e|),\\[1mm]
&&\partial_t (|\partial_t U^\e|)+\partial_x (F'(U^\e)|\partial_t U^\e|)
  \le \varepsilon \partial_{xx} (|\partial_t U^\e|)
\end{eqnarray*}
in the sense of distributions, leading to the BV--estimate.

Then the compactness theorem in BV implies the strong convergence of $U^\e(t,x)$.

\smallskip
Similar arguments can yield the $L^1$--equicontinuity of $U^\e$ directly,
which is also a corollary of the $L^1$-stability and the comparison
principle via Kruzhkov's method \cite{Kruzkov}.

\smallskip
The same arguments also work for multidimensional scalar conservation
laws ({\it cf.} \cite{Kruzkov,Volpert}); also see \cite{ChenChristoforou}
for scalar conservation laws with memory.

%=======================================================================

\subsubsection{Hyperbolic systems of conservation laws: $BV$--estimate via Glimm's approach}

Glimm \cite{Glimm}
first developed a random choice method, the Glimm scheme, and derived the $BV$-estimate
of the corresponding Glimm approximate solutions, based on
the Glimm functional and corresponding wave interaction estimates.
The techniques developed have been successfully employed to
establish the global existence of solutions in $BV$
and analyze the behavior of solutions in $BV$ (structure, uniqueness, stability,
and asymptotic behavior of solutions in $BV$) when the total variation of the initial data is small.
Also see Glimm-Lax \cite{Glimm-Lax}, DiPerna \cite{DiPerna77},
Liu \cite{Liu81},  Dafermos \cite{Dafermos},
and the references cited therein.

\smallskip
{\bf Theorem}: {\it For a strictly hyperbolic system \eqref{2.1}
on {$U$} in a neighborhood of
a compact set {$K\subset\R^N$}, there exist
constants $\delta>0$ and $C$ such that, if
\begin{equation}\label{2.3a}
{\rm Tot.Var.}\{U_0\}<\delta, \qquad \lim_{x\to -\infty}U_0(x)\in K,
\end{equation}
then there exists a global solution
$U(t,x)$ such that
$$
{\rm Tot.Var.}\{U(t,\cdot)\}\le C\, {\rm Tot.Var.} \{U_0\}.
$$
}

Glimm's approach has been further employed to handle the front-tacking method
and developed to analyze the $L^1$--stability of global solutions obtained by
either the Glimm scheme or the front tracking method.
See Bressan \cite{Bressan-book}, Dafermos \cite{Dafermos}, Holden-Risebro \cite{Holden-Risebro},
Liu-Yang \cite{LiuYang}, LeFloch \cite{LeFloch-book},
and the references cited therein.
The approach has also been developed to analyze the well-posedness for two-dimensional steady supersonic
Euler flows past a Lipschitz wedge in \cite{CZZ,ChenLi}.

\subsubsection{Hyperbolic system of conservation laws: BV--estimate for the artificial viscosity method}

Consider the following Cauchy problem for one-dimensional nonlinear hyperbolic systems
of conservation laws with vanishing artificial viscosity ({\it i.e.} $D(U)=I_{N\times N}$):
\begin{equation}\label{2.4.1}
\del_t U + \del_x F(U)
= \varepsilon \partial_{xx}U
\end{equation}
and the initial data: $U(0,x)=U_0(x)\in BV(\R^N)$.

\medskip
{\bf Theorem (Biachini-Bressan \cite{BB})}: {\it For a strictly hyperbolic system \eqref{2.1} on $U$
in a neighborhood of
a compact set $K\subset\R^N$, there exist
constants $\delta>0$ and $C_j, j=1,2,3$, such that, if $U_0$ satisfies \eqref{2.3a},
then, for any fixed $\e>0$, there exists a unique solution
{$U^\e(t,\cdot):=S_t^\e U_0(\cdot)$} of the Cauchy problem \eqref{2.4.1} such that

\begin{enumerate}\renewcommand{\theenumi}{\roman{enumi}}
\item[\rm (i)] {\it BV bound}:  $\, {\rm Tot.Var.}\{S_t^\e U_0\}\le C_1\, {\rm Tot.Var.}\{U_0\}$;

\item[\rm (ii)] {\it $L^1$--stability}:  $\, \|S_t^\e U_0-S_t^\e V_0\|_{L^1}\le C_2\, \|U_0-V_0\|_{L^1}$,

\medskip
$\qquad\qquad\quad\,\,\,$ $\|S_t^\e U_0-S_s^\e U_0\|_{L^1}\le C_3 \, \big(|t-s|+|\sqrt{\e t}-\sqrt{\e s}|\big)$.
\end{enumerate}
}

\smallskip
These imply the strong convergence and $L^1$-stability of the limit solution of \eqref{2.1}.

\medskip
The strategies to achieve the BV--estimate include the following steps:

\begin{enumerate}\renewcommand{\theenumi}{\roman{enumi}}
\item[(i)] Employ the {heat kernel} to estimate the solution for {$t\in [0, \tau_\e]$}:
$$
\|\partial_x U^\e(t,\cdot)\|_{L^1}\le \kappa \, \delta,
$$
where $\kappa$ is small, independent of $\e$ and $\delta$.

\item[(ii)] Decompose $\partial_x U^\e$ along a suitable basis of unit vectors
$\{\rr_1, \cdots, \rr_N\}$:
$$
\partial_x U^\e=\sum
v_i^\e \rr_i \quad
\mbox{(sum of gradients of viscous travelling waves)}.
$$

\item[(iii)] Derive a system of $N$ equations for these scalar components:
$$
\partial_t v_i^\e +\partial_x(\tilde{\lambda}_i v_i^\e)-\e\partial_{xx}v_i^\e=\phi_i^\e, \qquad
i=1, \cdots, N.
$$
Then, as the scalar case, we obtain that, for all $t\ge \tau_\e$,
$$
\|v_i^\e(t,\cdot)\|_{L^1}
\le \|v_i^\e(\tau_\e, \cdot)\|_{L^1}+\int_{\tau_\e}^\infty\int_{-\infty}^\infty|\phi_i^\e(t,x)|dxdt.
$$
\item[(iv)] Construct the basis $\{\rr_1, \cdots, \rr_N\}$ in an appropriate way so that, for $t\ge \tau_\e$,
$$
\int_{\tau_\e}^\infty\int_{-\infty}^\infty |\phi_i^\e(t,x)|dxdt\le \hat{C},
$$
which implies
$$
\mbox{\rm Tot.Var.}\{U^\e(t,\cdot)\}=\|U_x^\e(t,\cdot)\|_{L^1}\le \sum_i \|v_i^\e(t,\cdot)\|_{L^1}\le C,
$$
where $\hat{C}$ and $C$ are independent of $\e>0$.
\end{enumerate}

\smallskip
{\bf Remark 1.}
The results above still hold even for non-conservative strictly hyperbolic systems.
On the other hand, this approach requires both the artificial
viscosity ({\it i.e.} $D(U)=I_{N\times N}$)
and the total variation of the initial data sufficiently small.

\smallskip
{\bf Remark 2}. A longstanding open problem is the $BV$--estimate and convergence of vanishing
viscosity approximation $U^\e$ governed by the general form:
\begin{equation}\label{viscosity}
\del_t U^\varepsilon + \del_x F(U^\varepsilon)
= \varepsilon \partial_{x}\big(D(U^\e)\partial_x U^\varepsilon\big)
\end{equation}
for general viscosity matrices $D(U)$, including the Navier-Stokes viscosity matrices.
This especially includes the fundamental
problem in mathematical fluid dynamics, the inviscid limit of solutions of the Navier-Stokes equations
to the Euler equations for homentropic flow, via the BV--estimate, which is still open.

\subsection{Compactness and Convergence via Compensated Compactness}

We now discuss the compactness and convergence of exact/approximate solutions to conservation laws
via compensated compactness and related ideas,
which only require much weak bounds that may be obtained easily through
natural energy/entropy estimates as our examples below indicate.

\subsubsection{Scalar conservation laws}

Consider the Cauchy problem for scalar conservation laws \eqref{2.3} ($N=1$) with initial data:
$$
U|_{t=0}=U_0 \in L^\infty(\R).
$$
Then it can be easily shown that there exists $C$, independent of $\e$, such that
the viscous solutions $U^\eps$ satisfy the following natural estimates:

\begin{enumerate}
\renewcommand{\theenumi}{\roman{enumi}}
\item[(i)] {\it Maximum principle}: $\,\|U^\e\|_{L^\infty} \le C$  or $\|U^\e\|_{L^p}\le C$;

\smallskip
\item[(ii)] {\it Dissipation estimate}: $\, \|\sqrt{\e} U^\e_x\|_{L^2_{loc}}\le C$.
\end{enumerate}

The second estimate is a direct corollary of the natural energy estimate:
$$
\e |\partial_x U^\e|^2=-\partial_t\big(\frac{|U^\e|^2}{2}\big)-\partial_x\big(\int^{U^\e} wF'(w)dw\big)
+\e \partial_{xx}\big(\frac{|U^\e|^2}{2}\big).
$$

These estimates imply that, for any {$\eta\in C^2$} with entropy
flux $q(U)=\int^U\eta'(w)F'(w)dw$,
$$
\partial_t \eta(U^\e)+\partial_xq(U^\e) \quad\mbox{is compact in}\,\, H^{-1}_{loc}.
$$
Then the compensated compactness arguments yield
the weak continuity of $F(U^\e)$, or even the strong convergence of $U^\e(t,x)$ {\it a.e.}

For the convex case, Tartar \cite{Tartar} was the first to employ one entropy
pair $(\eta_*(U), q_*(U))=(U^2, 2\int^UwF'(w)dw)$ to conclude the strong convergence,
which initiated the successful applications of compensated compactness
to nonlinear hyperbolic conservation laws.
For the nonconvex case,
the entropy pair $(\eta_*(U), q_*(U))=(F(U), \int^U (F'(w))^2dw)$
also suffices to conclude the weak continuity with respect to
the general equation,
and the strong convergence when the equation is genuinely nonlinear for almost all $U$, as observed by
 Chen-Lu \cite{ChenLu} and Luc Tartar independently.
Also see
Schonbek \cite{Schonbek}, DiPerna \cite{DiPerna85a}, and Tadmor-Rascle-Bagnerini \cite{TRB}.

\medskip
The approach also applies to equation \eqref{viscosity} ($N=1$) with more general viscosity terms,
as well as scalar conservation laws with memory \cite{Dafermos88}.

\medskip
For these, the following Murat-Tartar's div-curl lemma plays an essential role:

\medskip
{\bf Div-Curl Lemma} (Tartar \cite{Tartar}, Murat \cite{Murat}):
{\it Let $\Omega\subset\R^d, d\ge 2,$ be open bounded. Let $p, q>1$ such
that {$\frac{1}{p}+\frac{1}{q}=1$}. Assume that, for any
$\varepsilon>0$, two vector fields
$$
{u^\varepsilon\in L^p(\Omega; \R^d)}, \qquad
{v^\varepsilon\in L^q(\Omega; \R^d)}
$$
satisfy the following:
\begin{enumerate}\renewcommand{\theenumi}{\roman{enumi}}
\item[\rm (i)]
\label{weakSol-def-i1}
{$u^\varepsilon\rightharpoonup u$ weakly in
$L^p(\Omega;\R^d)$ as $\varepsilon\to 0$;}
\item[\rm (ii)]
\label{weakSol-def-i2}
{$v^\varepsilon\rightharpoonup v$ weakly in
$L^q(\Omega;\R^d)$ as $\varepsilon\to 0$;}
\item[\rm (iii)]  {${\rm div}\, u^\varepsilon$ are confined in a
compact subset of $W^{-1, p}_{loc}(\Omega; \R)$;}
\item[\rm (iv)] {${\rm curl}\, v^\varepsilon$ are confined in a
compact subset of $W^{-1, q}_{loc}(\Omega; \R^{d\times d})$.}
\end{enumerate}
Then the scalar product of $u^\e$ and $v^\e$ are weakly continuous:
$$
u^\varepsilon\cdot v^\varepsilon\longrightarrow u\cdot v
$$
in the sense of distributions.}

\medskip
Various variations of this lemma for different applications/purposes
have been developed; see Tartar \cite{Tartar09}, Briane, Casado-Diaz and Murat \cite{BCM},
and the references cited therein.

\subsubsection{Hyperbolic systems of conservation laws: Compensated compactness
for the artificial viscosity method}

Consider system \eqref{2.4.1} with artificial viscosity.
Assume that there exists a strictly convex entropy function {$\eta_*(U), \nabla^2 \eta_*(U)>0$}.
In many cases, it can be shown that there exists {$C$} {independent of} {$\e$} such that

\begin{enumerate}\renewcommand{\theenumi}{\roman{enumi}}
\item[(i)]  {\it Invariant regions}:  $\quad \|U^\e\|_{L^\infty}\le C$;

\item[(ii)] {\it Dissipation estimate}:  $\quad \|\sqrt{\e}\partial_x U^\e\|_{L^2_{loc}}\le C$.
\end{enumerate}

In fact, the dissipation estimate is natural as the scalar case, directly from the energy estimate as follows:
$$
\e (\partial_x U^\e)^\top\nabla^2\eta_*(U^\e) \partial_x U^\e=-\partial_t\eta_*(U^\e)-\partial_x q_*(U^\e)
+\e \partial_{xx}\eta_*(U^\e).
$$
Then, for any {$\eta\in C^2$} with entropy flux {$q$}, {\it i.e.}, {$\nabla q(U)=\nabla\eta(U)\nabla F(U)$},
$$
\partial_t \eta(U^\e)+\partial_xq(U^\e) \quad\mbox{is compact in}\,\, H^{-1}_{loc}.
$$
The compensated compactness arguments can yield the strong convergence of $U^\e(t,x)$ when the system has strong
nonlinearity.

The similar compensated compactness arguments apply to the systems with more general viscosity matrices \eqref{viscosity}
for $\nabla^2\eta_*(U)D(U)\ge c_0>0$.
Another advantage of this approach is to allow the initial data of large oscillation without bounded variation.

\smallskip
In order to achieve the strong compactness, as first indicated by Tartar \cite{Tartar}, combining
the div-curl lemma (\cite{Tartar,Murat})
and the Young measure representation theorem ({\it cf.} Tartar \cite{Tartar}; also see \cite{Alberti-Muller,Ball2}),
we have the following commutation identity for the associated Young measure $\nu=\nu_{(t,x)}(\bla)$
(probability measure) for the sequence $U^\varepsilon(t,x)$:
\begin{equation}\label{commut-identity}
\begin{array}{ll}
&\langle\nu(\bla), \eta_1(\bla)q_2(\bla)-q_1(\bla)\eta_2(\bla)\rangle \\[1mm]
&\,=\, \langle\nu(\bla), \eta_1(\bla)\rangle\langle\nu(\bla),
 q_2(\bla)\rangle
 -\langle\nu(\bla), q_1(\bla)\rangle\langle\nu(\bla),
\eta_{2}(\bla)\rangle
\end{array}
\end{equation}
for any entropy pairs $(\eta_j, q_j), j=1,2$, and
$$
\partial_t\langle\nu, \eta(\bla)\rangle +\partial_x \langle\nu, q(\bla)\rangle\le 0
$$
in the sense of distributions
for any convex entropy pair $(\eta, q), \nabla^2\eta\ge 0$.
Then the main mathematical issue is whether $\nu$ is a Dirac measure.
The key point is the imbalance of regularity of the operator in the commutation
identity: The operator on the left is more regular than the one on the right
due to cancellation when the system has strong nonlinearity.
If so, the compactness of $U^\varepsilon(t,x)$ in $L^1$ follows.

\smallskip
For strict hyperbolicity with $N=2$, there are two families of entropy pairs determined by two arbitrary functions,
which yield an affirmative answer to the issue; see
DiPerna \cite{DiPerna83b}, Dafermos \cite{Dafermos}, Serre \cite{Serre87},
Morawetz \cite{Morawetz91}, Perthame-Tzavaras \cite{Perthame-Tzavaras}, and Chen-Li-Li \cite{ChenLiLi}.

\medskip
Further challenges include
nonstrictly hyperbolic systems,
viscosity matrices with $\nabla^2\eta_*(U)D(U)\ge 0$ but not positive definite,
and initial data of large oscillation with only energy bounds (without bounded variation or
$L^\infty$--uniform bound).
We now start with a fundamental example of nonstrictly hyperbolic systems.

\subsubsection{Homentropic Euler equations: Compensated compactness for the artificial viscosity method}

The homentropic Euler equations take the following form:
\begin{equation}\label{2.2.3.1}
\left\{\begin{array}{ll}
\del_t\rho + \del_x m =0, \\[1mm]
\del_t m + \del_x\big(\frac{m^2}{\rho}+p(\rho)\big)= 0,
\end{array}
\right.
\end{equation}
where $\rho$ is the density, $u=\frac{m}{\rho}$ is the fluid velocity that is well-defined when $\rho>0$,
$p=p(\rho)=\rho^2e'(\rho)$ is the pressure with internal energy $e(\rho)$,

For a polytropic perfect gas,
\begin{equation}\label{polytropic}
p(\rho)=\kappa\rho^\gamma, \,\,\,
e(\rho)=\frac{\kappa}{\gamma-1}\rho^{\gamma-1},
\end{equation}
where $\gamma>1$ is the adiabatic exponent,
and constant $\kappa$ in the pressure-density relation may be chosen as
$\kappa=\frac{(\gamma-1)^2}{4\gamma}$ without loss of generality.

\smallskip
One of the main difficulties for solving \eqref{2.2.3.1} is that
strict hyperbolicity fails when {$\rho\to 0$}.

\smallskip
An entropy function $\eta(\rho, m)$ is called a weak entropy if
$\eta(\rho, m)|_{\rho=0}=0$.
The weak entropy pairs can be represented as
\begin{equation}\label{WeakEntropy}
\eta^\psi(\rho,\rho u)
=\int_{\mathbb{R}}\chi(s)\psi(s)\,ds,\,\,\,
q^\psi(\rho,\rho u)
=\int_{\mathbb{R}}\big(\theta s{}+{}(1-\theta)u\big)\chi(s)\psi(s)\,ds
\end{equation}
for any $C^2$-test function $\psi(s)$, where $\chi(s)$ is the weak entropy kernel:
\begin{equation}\label{kernel}
\chi(s):=[\rho^{2\theta}{}-{}(u-s)^2\big]_+^\lambda,  \qquad
\theta=\frac{\gamma-1}{2}, \lambda=\frac{3-\gamma}{2(\gamma-1)}.
\end{equation}

The mechanical energy--energy flux pair
$(\eta_*, q_*)$:
$$
\eta_*(\rho,m)=\frac{1}{2}\frac{m^2}{\rho}+ \rho e(\rho), \qquad
q_*(\rho,m)=\frac{1}{2}\frac{m^3}{\rho^2}+m \big(e(\rho)+\frac{p}{\rho}\big)
$$
is a convex entropy pair for \eqref{2.2.3.1}.

\smallskip
Consider the homentropic Euler equations with artificial viscosity:
\begin{equation}\label{2.2.3.2}
\left\{\begin{array}{ll}
\del_t\rho + \del_x m =\e \del_x^2\rho, \\[1mm]
\del_t m + \del_x\big(\frac{m^2}{\rho}+p(\rho)\big)=\e\del_x^2 m.
\end{array}
\right.
\end{equation}
It can be shown for system \eqref{2.2.3.2} that
there exists $C>0$, independent of $\e>0$, such that
\begin{enumerate}\renewcommand{\theenumi}{\roman{enumi}}
\item[(i)] {\it Invariant regions for the $L^\infty$--estimate}:
$$
0\le \rho^\e(t,x)\le C, \qquad
|m^\e(t,x)|\le C \rho^\e(t,x) \qquad\quad\mbox{\it a.e.}
$$
\item[(ii)] {\it Dissipation estimate}:
$$
\sqrt{\e}\|\del_x(\rho^\e, m^\e)\|_{L^2([0,T]\times\R)}\le C,
$$
via the mechanical energy pair
$(\eta_*, q_*)$ that is strictly convex for $1<\gamma\le 2$,
and convex for $\gamma>2$ for which a corresponding
weighted dissipation estimate can be obtained.
\end{enumerate}

These estimates yield that, for any $C^2$ weak entropy
pair $(\eta, q)$,
$$
\del_t\eta(\rho^\e, m^\e)+\del_xq(\rho^\e, m^\e)
\quad\mbox{is compact in $H^{-1}_{loc}$}.
$$
Then the convergence problem for $(\rho^\e, m^\e)$ is reduced to
the reduction problem for a measure-valued solution $\nu_{t,x}$:

\medskip
{\it If $\mbox{supp} \, \nu_{t,x}$ is bounded,
then
\begin{equation}\label{2.11a}
\nu_{t,x}=\nu_{(\rho(t,x), m(t,x))},
\end{equation}
that is, $\,\, (\rho^\e(t,x), m^\e(t,x))\,\,\to \, (\rho(t,x),
m(t,x))$ a.e. $(t,x)$.
}

\medskip
This problem has been solved by
DiPerna \cite{DiPerna83a} for $\gamma=\frac{N+2}{N}$, $N\ge 5$ odd,
Ding-Chen-Luo \cite{DCL} and Chen \cite{Chen3} for $\gamma \in
(1, \frac{5}{3}]$,
Lions-Perthame-Tadmor \cite{LPT} for $\gamma\ge3$,
Lions-Perthame-Souganidis \cite{LPS} for $\gamma \in (\frac{5}{3}, 3)$,
and
Chen-LeFloch \cite{CLF} for general pressure laws.
The key point is to employ effectively the weak entropy pairs
in the commutation identity \eqref{commut-identity} for the associated Young measure
$\nu_{t,x}$ with compact support.

\smallskip
The convergence of related numerical methods with corresponding numerical viscosity
matrices including the Lax-Friedrichs scheme and Godunov scheme has also been established
in Ding-Chen-Luo \cite{DCL}; also see Chen \cite{Chen05}.

\smallskip
The isothermal case $\gamma=1$ has also been handled by Huang-Wang \cite{Huang-Wang};
also see LeFloch-Shelukhin \cite{LeFloch-Shelukhin}.

\smallskip
Some further important problems include
the inviscid limit from the compressible Navier-Stokes equations to the compressible Euler equations
(see \S 2.3)
and the existence of global spherically symmetric solutions to the compressible Euler equations
(see \S 2.4).

\subsection{Navier-Stokes Equations: Inviscid Limit}

Consider the Cauchy problem:
\begin{equation}\label{2.13a}
\left\{
\begin{array}{ll}
\rho_t{}+{}(\rho u)_x{}={}0,\\[2mm]
(\rho u)_t{}+{}(\rho u^2{}+{}p)_x{}={}\e u_{xx}{},\\
\end{array}
\right.
\end{equation}
with the initial conditions:
\begin{equation}\label{2.14a}
(\rho, u)|_{t=0}=(\rho_0^\e(x), u_0^\e(x)),\quad\,
\lim_{x\to\pm\infty}(\rho_0^\e(x), u_0^\e(x))=(\rho^\pm, u^\pm),
\end{equation}
where
$(\rho^\pm, u^\pm)$ are constant end-states with $\rho^\pm>0$,
and the viscosity coefficient $\e\in (0,\e_0]$ for some fixed $\e_0$.

\smallskip
The existence of $C^2$--solutions $(\rho^\e, u^\e)(t,x)$ for large initial
data was obtained by Kanel \cite{Kanel} for the same ending states
and by Hoff \cite{Hoff} for different ending states.

\medskip
{\bf Inviscid Limit Problem:} {\it Does the solution
sequence $(\rho^\e, u^\e)(t,x)$
of system \eqref{2.13a} strongly converge
to a solution to the homentropic
Euler equations \eqref{2.2.3.1}
when $\e\to 0$?}

\medskip
This problem has been addressed by Gilbarg \cite{Gilbarg}, Hoff-Liu \cite{Hoff-Liu},
and G\`{u}es-M\'{e}tivier-Williams-Zumbrun \cite{GMWZ} for
some physical cases with special structure for which the limit solution
contains only one shock.

\smallskip
For the general case, several new difficulties arise, which include
\begin{enumerate}\renewcommand{\theenumi}{\roman{enumi}}
\item[(i)] No invariant regions:   Only energy norms;
\item[(ii)] Direct derivative estimates only partially:
$\|\sqrt{\e}\del_x u^\e\|_{L^2_{loc}}\le C$;

\item[(iii)] No {\it a priori} bounded support of the measure-valued solution
$\nu_{t,x}$.
\end{enumerate}

Nevertheless, the following theorem has been established.

\medskip
{\bf Theorem} (Chen-Perepelitsa \cite{CP}):
{\it Let the initial functions $(\rho_0^\e, u_0^\e)$ satisfy
\begin{eqnarray*}
&&\int_{-\infty}^\infty \Phi_*(\rho_0^\e(x), m^\e_0(x))
dx\le E_0<\infty,\\
&&\int_{-\infty}^\infty\Big(\e^2\frac{|\rho_{0,x}^\e(x)|^2}{\rho_0^\e(x)^3}{}+{}2\e
\frac{|\rho_{0,x}^\e(x)u_0^\e(x)|}{\rho_0^\e(x)}+ \rho_0^\e(x)|u_0^\e(x)| \Big)dx \le E_1<\infty,\\
\end{eqnarray*}
where
$\, \Phi_*(\rho,m)=\eta_*(\rho,m)-\eta_*(\bar{\rho}, \bar{m})-\nabla\eta_*(\bar{\rho}, \bar{m})\cdot (\rho-\bar{\rho}, m-\bar{m})\ge 0$
for
$\bar{m}=\bar{\rho}\bar{u}$,
$(\bar{\rho}, \bar{u})$ is a pair of smooth monotone functions satisfying
$(\bar{\rho}(x), \bar{u}(x))=(\rho^\pm, u^\pm)$ when $\pm x\ge L_0$ for some large $L_0 > 0$,
and both $E_0$ and $E_1$ are {independent of} {$\e$}.
Let $(\rho^\e, m^\e), m^\e=\rho^\e u^\e$, be the solution of
the Cauchy problem for the Navier-Stokes equations \eqref{2.13a}
for each fixed $\e>0$.
Then, when {$\e\to 0$}, there exists a subsequence
of $(\rho^\e, m^\e)$ that
converges strongly almost everywhere to a finite-energy
solution $(\rho, m)$ to the Cauchy problem for
the homentropic Euler equations
\eqref{2.2.3.1}
for any $\gamma>1$.
}

\medskip
The strategies for this include the following steps.

\begin{enumerate}\renewcommand{\theenumi}{\roman{enumi}}
\item[(i)] Derive the finite-energy bound and higher integrability bound (replacing $L^\infty$ bound);
\item[(ii)] Derive new derivative estimate for $\e\del_x\rho^\e$;
\item[(iii)] Show the $H^{-1}$--compactness of weak entropy dissipation measures only
for weak entropy pairs with compactly supported $C^2$--test
functions;
\item[(iv)] Prove that any connected component of support of the
measure-valued solution $\nu_{t,x}$ must be bounded,
which reduces to the case when the support of $\nu_{t,x}$ is bounded as in \S 2.2.3.
\end{enumerate}

\smallskip
To achieve these, the following key estimates of solutions to the Navier-Stokes equations
are essential:
There exist $C_1>0$
and $C_2=C_2(E_0, E_1, K, \gamma, t)$ independent of $\e$
for any compact set {$K\subset \R$} such that, for any
$t>0$,

\begin{enumerate}\renewcommand{\theenumi}{\roman{enumi}}

\item[(i)] Energy estimate:
$$
\int_{-\infty}^\infty \Phi_*(\rho^\e(t,x), m^\e(t,x))\, dx
+{}\int_0^t\int_{-\infty}^\infty\, \e|u_x^\e|^2\, dxd\tau {}\leq{} E_0;
$$

\item[(ii)] New derivative estimate for the density:
$$
\e^2\int\frac{|\rho_x^\e(t,x)|^2}{\rho^\e(t,x)^3}dx{}
+{}\e\int_0^t\int_{-\infty}^\infty (\rho^\e)^{\gamma-3}|\rho_x^\e|^2\, dxd\tau
{}\leq{} C_1(E_0+E_1).
$$

\item[(iii)] Higher integrability bound:
\begin{equation*}
\int_0^t\int_K
\big(\rho^\e|u^\e|^3{}+{}(\rho^\e)^{\gamma+\theta}
     +(\rho^\e\big)^{\gamma+1}\big)\, dxd\tau{}\leq{} C_2.
\end{equation*}
\end{enumerate}

The higher integrability estimate (iii) is motivated by the
related work by Lions-Perthame-Tadmor \cite{LPT} and LeFloch-Westdickenberg \cite{LW}.
For some related earlier work on the convergence of approximate solutions in the $L^p$--framework,
see Serre-Shearer \cite{SS} for a $2\times 2$ system of
elasticity with severe growth conditions,
and LeFloch-Westdickenberg \cite{LW} for the convergence of
approximate solutions with full dissipation in the energy norms for the
homentropic Euler equations with $\gamma\in (1, \frac{5}{3}]$.

\medskip
Let $\nu_{t,x}$ be the Young measure determined by the
solutions of the Navier-Stokes equations \eqref{2.13a}.
Then $\nu_{t,x}$
is confined
by
\begin{eqnarray*}\label{commute}
&&
\theta(s_2-s_1)\big(\overline{\chi(s_1)\chi(s_2)}-
\overline{\chi(s_1)}\,\, \overline{\chi(s_2)}\big)
\\&&\,\,
=(1-\theta)\big(\overline{u\chi(s_2)}\,\, \overline{\chi(s_1)}
-\overline{u\chi(s_1)}\, \,\overline{\chi(s_2)}\big) \quad\mbox{for
{\it a.e.} $s_1, s_2\in \R$,}
\end{eqnarray*}
for the entropy kernel $\chi(s):=[\rho^{2\theta}{}-{}(u-s)^2\big]_+^\lambda$
with $\theta=\frac{\gamma-1}{2}$ and $\lambda=\frac{3-\gamma}{2(\gamma-1)}$,
where $\overline{f(s)}:=\langle \nu_{t,x}, f(s; \rho, u)\rangle$.

\medskip
The goal is to establish that the Young measure is a Dirac mass in the phase plane
for $(\rho, m)$.
The new difficulty is now that supp $\nu_{t,x}$ {\it is unbounded} in general.

\medskip
We divide the proof into three cases:
\medskip

%%%%%%%%%%%%%%%%%%%%%%%%%%%%%%%%%%
{\bf Case 1}: $\gamma=3$. The same argument for the bounded support of $\nu_{t,x}$ applies as in \cite{LPT}.
In this case, $\theta=1$ and the commutation relation
becomes
\[
\overline{\chi(s_1)\chi(s_2)}{}={}\overline{\chi(s_1)}\,\,\overline{\chi(s_2)},
\]
which implies
$\overline{\chi(s)}^2{}={}\overline{\chi(s)^2}$
by taking $s_1=s_2$. That is,
$$
\langle \nu_{t,x}, \big(\chi(s)-\overline{\chi(s)}\big)^2\rangle=0
\qquad\mbox{for any $s\in \R$}.
$$
This implies that
$\nu$ must be a Dirac mass on the set $\{\rho>0\}$ or be
supported completely in the vacuum $V=\{\rho=0\}$, that is,
the measure-valued solution $\nu_{t,x}$ is a Dirac mass
\eqref{2.11a}
in the
phase plane for $(\rho, m)$.

\medskip
{\bf Case 2}: {$\gamma>3$}.
Let $A:={}\cup \{ (u-\rho^\theta,\rho^\theta+u)\,:\,(\rho,u){}\in{}
\mbox{supp}\, \nu \}$.
Let $J=(s_-, s_+)$ be any connected component
of $A$.
Note that $\mbox{supp}\, \chi(s){}={}\{ (\rho,u)\,:\,
u-\rho^\theta\le s\le u+\rho^\theta\}.$

\smallskip
{\bf Claim}: $\,$ {\it The connected component $J$ is bounded for $\gamma>3$}.

\smallskip
On the contrary, let $\inf\{s{}:{}s\in J\}=-\infty$.
Our strategy is to fix $M_0$ first such that $M_0+1\in J$ and restrict $s_2\in (M_0, M_0+1)$,
and then choose sufficiently small $s_1\le -2|M_0|$ to reach the contradiction.

\smallskip
To achieve this, two following estimates are essential:

\medskip
(i) $\int_{M_0}^{M_0+1}\frac{\overline{\chi(s_1)\chi(s_2)}}{\overline{\chi(s_1)}}ds_2
\le C(\lambda) |s_1|^\lambda$ for $\lambda<0$, which is our key new observation.

\medskip
(ii) $\frac{\overline{\chi(s_1)\chi(s_2)}}{\overline{\chi(s_1)}}{}\geq{}\overline{\chi(s_2)}$
{\it a.e.}  $s_1,s_2\in J,s_1<s_2$,
by employing Lions-Perthame-Tadmor's argument in \cite{LPT}.

\smallskip
Combining the two estimates, we have
$$
\int_{M_0}^{M_0+1}\frac{\overline{\chi(s_1)\chi(s_2)}}{\overline{\chi(s_1)}}ds_2
\ge \int_{M_0}^{M_0+1}\overline{\chi(s_2)}ds_2 =C(M_0,\lambda)>0,
$$
which implies that, when $s_1\to-\infty$,
$$
0<C(M_0,\lambda)=\int_{M_0}^{M_0+1}\frac{\overline{\chi(s_1)\chi(s_2)}}{\overline{\chi(s_1)}}ds_2
\le C(\lambda) |s_1|^\lambda  \to 0.
$$
This arrives at the contradiction.

\smallskip
The case when $J$ is unbounded from above can be treated similarly.

\smallskip
This indicates that any connected component
$J$ of the support of the Young
measure $\nu$ is bounded for $\gamma>3$,
which reduces to the Lions-Perthame-Tadmor's case for $\gamma>3$ in \cite{LPT}.

\medskip
{\bf Case 3}: $\gamma\in (1,3)$.
On the contrary, suppose that a connected component $J$ of the support is unbounded from below.

Let $M_0{}={}\sup\{s{}:{}s\in J\}\in(-\infty,\infty]$.
Let $\xi_1, \xi_2, \xi_3\in (-\infty, M_0)$ with
$\xi_1<\xi_2<\xi_3$.
The commutation relation leads to
\begin{eqnarray}
&&(\xi_2-\xi_1)\frac{\overline{\chi(\xi_1)\chi(\xi_2)}}{\overline{\chi(\xi_1)}}
{}+{}(\xi_3-\xi_2)\frac{\overline{\chi(\xi_3)\chi(\xi_2)}}{\overline{\chi(\xi_3)}}
\nonumber\\
&& =(\xi_3-\xi_1)\overline{\chi(\xi_2)}\frac{\overline{\chi(\xi_1)\chi(\xi_3)}}
{\overline{\chi(\xi_1)}\,\overline{\chi(\xi_3)}}. \label{commut-1}
\end{eqnarray}
Differentiating this equation in $\xi_2$ and dividing by
$(\xi_3-\xi_1)$, we obtain
\begin{eqnarray}
\overline{\chi'(\xi_2)}\frac{\overline{\chi(\xi_1)\chi(\xi_3)}}
{\overline{\chi(\xi_1)}\,\overline{\chi(\xi_3)}}
&=&\frac{\xi_2-\xi_1}{\xi_3-\xi_1}\frac{\overline{\chi(\xi_1)\chi'(\xi_2)}}{\overline{\chi(\xi_1)}}
{}+{}\frac{\xi_3-\xi_2}{\xi_3-\xi_1}\frac{\overline{\chi(\xi_3)\chi'(\xi_2)}}{\overline{\chi(\xi_3)}}\nonumber\\
&&+\frac{1}{\xi_3-\xi_1}\frac{\overline{\chi(\xi_1)\chi(\xi_2)}}{\overline{\chi(\xi_1)}}
-{}\frac{1}{\xi_3-\xi_1}\frac{\overline{\chi(\xi_3)\chi(\xi_2)}}{\overline{\chi(\xi_3)}}.\nonumber\\
\label{commut-2}
\end{eqnarray}

Our strategy is to take $\xi_1\to-\infty$ first and show then that the
left-hand side has a smaller order than the right-hand side
to arrive at the contradiction.

\medskip
To do this, we divide the argument into five steps:

\smallskip
(i).  Show the estimate:
$$
\frac{\overline{\chi(\xi_1)\chi(\xi_3)}}{\overline{\chi(\xi_1)}
\,\,\overline{\chi(\xi_3)}}{}\geq{}1 \qquad \mbox{for any $s_1,s_3\in J$},
$$
by employing Lions-Perthame-Tadmor's argument in \cite{LPT}.

\medskip
(ii). Show that $\overline{\chi(\xi)}{}\geq{}0$, but is not identically zero, and
$\overline{\chi(\xi)}{}\to{}0\,\,$  as $\,\, \xi\to\inf J,\,\sup
J$.
This yields that there exists $\xi_2$ such that
$\overline{\chi'(\xi_2)}{}>{}0$ and $\overline{\chi(\xi_2)}{}>{}0$.

\medskip
(iii). Let $\xi_3>\xi_2$ be the points such that
$\overline{\chi(\xi_3)}>0$, and let $\xi_1\to-\infty$.
From the first identity \eqref{commut-1},
$$
\frac{\overline{\chi(\xi_1)\chi(\xi_2)}}{\overline{\chi(\xi_1)}}
{}={} \overline{\chi(\xi_2)}
\frac{\overline{\chi(\xi_1)\chi(\xi_3)}}
{\overline{\chi(\xi_1)}\,\,\overline{\chi(\xi_3)}}{}+{}o(1)
\qquad \mbox{as $s_1\to-\infty$}.
$$

(iv). Show that $[\chi'(s)]_+\le \frac{2\lambda}{s-s_1}\chi(s)$.

\medskip
(v). From the second equation \eqref{commut-2},
by throwing away the negative
terms, we have
\begin{eqnarray*}
\overline{\chi'(\xi_2)}
\frac{\overline{\chi(\xi_1)\chi(\xi_3)}}{\overline{\chi(\xi_1)}\,\,\overline{\chi(\xi_3)}}{}
\leq \frac{2\lambda+1}{\xi_3-\xi_1}\frac{\overline{\chi(\xi_1)\chi(\xi_2)}}{\overline{\chi(\xi_1)}}
{}+{}o(1),
\end{eqnarray*}
which implies
$$
\Big(
\overline{\chi'(\xi_2)}{}-{}\frac{2\lambda+1}{s_3-s_1}\overline{\chi(s_2)}\Big)
\frac{\overline{\chi(\xi_1)\chi(\xi_3)}}{\overline{\chi(\xi_1)}\,\,\overline{\chi(\xi_3)}}{}
\leq{}o(1).
$$
This arrives at the contradiction as $s_1\to -\infty$.

\medskip
Another different proof is given by LeFloch-Westdickenberg \cite{LW} for $1<\gamma\le \frac{5}{3}$.
The inviscid limit of the viscous shallow water equations to the Saint-Venant system has also been
established in Chen-Perepelitsa \cite{CP2}.

\subsection{Spherically Symmetric Solutions to the Multidimensional Homentropic Euler Equations}

The homentropic Euler equations for multidimensional compressible fluids take the following form:
\begin{equation}\label{ss-1}
\begin{cases}
\rho_t +\nabla_\x(\rho \vv)=0,\\[1mm]
(\rho \vv)_t+\nabla_\x(\rho \vv\otimes\vv)+ \nabla_\x p=0,
\end{cases}
\end{equation}
where $\x=(x_1, \dots, x_d)\in \R^d$, $\nabla_\x$ is the gradient with respect to $\x\in\R^d$,
and $\vv=(v_1, \dots, v_d)\in \R^d$ is the velocity.
The pressure-density constitutive relation (by scaling) satisfies \eqref{polytropic}.

\smallskip
We seek the spherically symmetric solutions with form:
\begin{equation}\label{ss-2}
\rho(t,\x)=\rho(t,r), \qquad \vv(t,\x)=u(t,r)\frac{\x}{r}, \qquad\quad r=|\x|.
\end{equation}
Then the functions $(\rho, m)=(\rho, \rho u)$ are governed by
\begin{equation}\label{ss-3}
\left\{ \begin{array}{l}
\rho_t{}+{}m_r{}+{}\frac{d-1}{r}m{}=0,\\[1mm]
m_t{}+{}(\frac{m^2}{\rho}{}+{}p(\rho))_r{}+{}\frac{d-1}{r}\frac{m^2}{\rho}
{}=0.
\end{array}
\right.
\end{equation}

For the defocusing case,
the existence of expanding spherically symmetric solutions with the following bounds:
$$
0\le  \rho(t,r)^{\frac{\gamma-1}{2}}\le u(t,r)\le C<\infty,
$$
has been constructed, provided that the initial functions have the same bounds,
in Chen \cite{Chen97}.

\smallskip
For the focusing case, the singularity of imploding self-similar spherically symmetric
solutions
has been discussed in
\cite{CFr,Guderley,Ro,Wh}.
It is indicated indeed in Rauch \cite{Rauch} that
there is no $BV$ or $L^\infty$ bound for the imploding solutions in general.

\smallskip
A longstanding open problem is whether the concentration phenomenon occurs at the origin, that is,
whether the density $\rho$ develops a measure at the origin.
In Chen-Perepelitsa \cite{CP3},
we have developed a method of vanishing artificial viscosity
to prove that the vanishing viscosity limit solution does not form concentration
at the origin, but has a bounded total energy.
More precisely, we construct a sequence of vanishing viscosity solutions to the
following initial-boundary problem:
\begin{equation}\label{ss-4} \left\{
\begin{array}{l}
\rho_t{}+{} m_r{}+\frac{d-1}{r}m={}\e \big(\rho_{rr}+\frac{d-1}{r}\rho_r\big),\\[1mm]
m_t{}+{}\big(\frac{m^2}{\rho}{}+{}p_\delta(\rho)\big)_r{}+\frac{d-1}{r}\frac{m^2}{\rho}={}\e \big( m_{rr}+\frac{d-1}{r}m\big)_r,
\end{array}
\right.
\end{equation}
with appropriate approximate initial data:
\begin{equation}\label{ss-5}
(\rho, m)|_{t=0}=(\rho_0^\e(r),m_0^\e(r))\to (\rho_0, m_0) \qquad a.e. \,\,\,\,\mbox{as } \, \e\to 0,
\end{equation}
and the boundary condition:
\begin{equation}\label{ss-6}
(\rho_r, m)|_{r=a(\e)}=(0, 0), \quad (\rho, m)|_{r=b(\e)}=(\bar{\rho}(\e), 0),
\end{equation}
where $(\rho_0,m_0)$ is the initial data for the spherical symmetric solution to system \eqref{ss-3},
$p_\delta(\rho)=\delta\rho^2+\kappa \rho^\gamma$ with $\delta=\delta(\e)$, and
$a(\e), b(\e), \bar{\rho}(\e)$, and $\delta(\e)$ are positive with
$a(\e)\to 0$, $b(\e)\to \infty$, $(\bar{\rho}(\e), \delta(\e))\to (0,0)$,
as well as certain combinations of $(b(\e), \bar{\rho}(\e), \delta(\e))$ tending to $0$,  as $\e\to 0$
({\it cf.} \cite{CP3}).
Then we have

\medskip
{\bf Theorem} (Chen-Perepelitsa \cite{CP3}).
{\it Let the initial functions $(\rho_0, m_0)$ for system \eqref{ss-3} satisfy the finite-energy
conditions. Then
\begin{enumerate}\renewcommand{\theenumi}{\roman{enumi}}
\item[\rm (i)] For sufficiently small fixed $\e>0$, there exists a global viscous solution
 $(\rho^\e, m^\e)$ to the initial-boundary value problem \eqref{ss-4}--\eqref{ss-6}
 satisfying that, for any compact set $K\subset \R_+$ and $T>0$,
there exists $C_T>0$ independent of $\e>0$ such that, for any $0<t\le T$,
\begin{eqnarray*}
&&\int_{a(\e)}^{b(\e)}\big(\frac{1}{2}\rho^\e (u^\e)^2+\rho^\e e(\rho^\e)\big)(t,r)\,r^{d-1}dr\\
&&\quad +\e\int_0^t\int_{a(\e)}^{b(\e)}\big((\rho^\e)^{\gamma-2}|\rho^\e_r|^2+\rho^\e|u^\e_r|^2+\frac{\rho^\e (u^\e)^2}{2r^2}\big)r^{d-1}drdt
\le C_T,
\end{eqnarray*}
and
$$
\int_0^t\int_K
\big(\rho^\e|u^\e|^3{}+{}(\rho^\e)^{\gamma+\theta} +(\rho^\e)^{\gamma+1}\big)\,r^{d-1} dr d\tau{}\leq{} C_T;
$$

\item[\rm (ii)] When $\e\to 0$, there exists a subsequence
of $(\rho^\e, m^\e)$ that converges
strongly almost everywhere to a finite-energy spherically symmetric
solution $(\rho, m)$ to system \eqref{ss-3} for any $\gamma>1$ with initial data $(\rho_0, m_0)$.
\end{enumerate}
}

The key ingredients are the uniform {\it a priori} estimates in (i) in combination with the reduction of the
corresponding Young measure discussed in \S 2.3.

\smallskip
Recently, we have also solved and/or made progress on several fundamental problems in
nonlinear partial differential equations by employing the
viscosity method. These include
vanishing viscosity approximation for transonic flow in Chen-Slemrod-Wang \cite{CSW08}
(also see Morawetz \cite{Mor85,Mor95}),
and subsonic-sonic limit of exact/approximate solutions to
the full Euler equations for multidimensional steady compressible fluids
in Chen-Huang-Wang \cite{CHW14}.

\section{Weak Continuity and Rigidity of the Gauss-Codazzi-Ricci System and Corresponding Isometric Embeddings}

The isometric embedding problem is a longstanding fundamental problem in
differential geometry.
As is well-known from differential geometry, given a surface, we can compute its
metric $\{g_{ij}\}$ and associated first fundamental form:
$$
I=\sum g_{ij}dx^i dx^j,
$$
and its curvatures determined by the second fundamental form:
$$
II=\sum h_{ij} dx^i dx^j.
$$
Then a natural mathematical question is as follows:

\medskip
{\bf Isometric Embedding Problem}: {\it Given a metric $\{g_{ij}\}$,
can we find a surface in the Euclidean space with
the given metric $\{g_{ij}\}$?}

\medskip
In other words, we seek a map
$\mathbf{r}: \Omega\to \mathbb{R}^N$ such
that
$$
d\mathbf{r}\cdot d\mathbf{r}=\sum_{i,j=1}^N g_{ij}dx^i dx^j
$$
in the local coordinates,
that is, $\d_{x^i}{\bf r}\cdot\d_{x^j}{\bf r}=g_{ij}$
   so that $(\d_{x^i}{\bf r}, \d_{x^j}{\bf r}), i\ne j$, in
$\mathbb{R}^d$ are linearly independent.

This is an inverse problem, which is a realization question for given
an abstract metric $\{g_{ij}\}$.
A further question is whether we can produce even more sophisticated surfaces
or thin
sheets for applications.
These questions are truly fundamental, not only in mathematics such as differential geometry and topology,
but also in many applications such as the understanding evolution of sophisticated shapes of surfaces
or thin sheets in nature including elastic materials, protein folding in biology and algorithmic origami,
as well as design, visual arts, among others.

\smallskip
The mathematical study of this problem has a long history, including the early important work by
Schlaefli (1873), Darboux (1894), Hilbert (1901), Weyl (1916), Janet (1926-27),
Cartan (1926-27); also see Han-Hong \cite{HanHong} and the references cited therein.
In particular, Nash \cite{Nash56} established the Nash isometric embedding theorem
(also called $C^k$--embedding theorem, $k\ge 3$):

\medskip
{\it Every $n$-dimensional Riemannian manifold (analytic or $C^k$, $k\ge 3$) can be $C^k$--isometrically imbedded in the
Euclidean space $\R^d$
with $d=2s_n+4n$ for the compact case and $d=(n+1)(3s_n+4n)$ for the noncompact case, where
$s_n=\frac{n(n+1)}{2}$ is the Jenet dimension ({\it cf.} {\rm \cite{Janet}}).}

\medskip
The results
were further improved with lowerer target dimensions by Gromov \cite{Gromov86}
with $d=s_n+2n+3$ and G\"{u}nther \cite{Gunther89} with $d=\max\{s_n+2n, s_n+n+5\}$.

\medskip
The following further problems are important for applications:
\begin{enumerate}\renewcommand{\theenumi}{\roman{enumi}}
\item[(i)] {\it Rigidity of isometric embeddings}: Is a weak limit of a sequence of isometric embeddings in some topology still an isometric
embedding?

\item[(ii)] {\it Lowerest target dimension for global isometric embeddings}, which is expected to be the Janet dimension $d=s_n$;

\item[(iii)] {\it Optimal or assigned regularity such as $C^{1,1}$, $W^{2, p}$, and $BV^1$}.
The regularity issue is quite sensitive. For example, Efimov's example in \cite{Efimov-a}
indicates that there is no $C^2$--isometric embedding
when $n=2$ and $d=s_n=3$.
\end{enumerate}

For $n=2$ and $d=3$, the fundamental theorem in differential geometry indicates that

\medskip
{\it There exists a surface in $\R^3$ whose first and second fundamental
forms are $I$ and $I\!I$, if the coefficients $\{g_{ij}\}$ and
 $\{h_{ij}\}$ of the two given quadratic forms $I$ and $I\!I$, $I$
being positive definite, satisfy the Gauss-Codazzi system.
That is, given $\{g_{ij}\}$, the second fundamental coefficients $\{h_{ij}\}$
are determined by the Codazzi equations (compatibility):
\begin{equation} \label{g1}
\left\{\begin{aligned}
&\d_x{M}-\d_y{L}=L\G^{(2)}_{22}-2M\G^{(2)}_{12}+N\G^{(2)}_{11}, \\
&\d_x{N}-\d_y{M}=-L\G^{(1)}_{22}+2M\G^{(1)}_{12}-N\G^{(1)}_{11},
\end{aligned}\right.
\end{equation}
subject to the Gauss equation ({\it i.e.}, the Monge-Amp\`{e}re type
constraint):
\begin{equation}\label{g2}
LN-M^2=K,
\end{equation}
where
$$
L=\frac{h_{11}}{\sqrt{|g|}}, \,\,
  M=\frac{h_{12}}{\sqrt{|g|}}, \,\,
  N=\frac{h_{22}}{\sqrt{|g|}}, \,\, |g|=g_{11}g_{22}-g_{12}^2,
$$
$\G^{(k)}_{ij}$ are the Christoffel symbols, depending
on $g_{ij}$ up to their first derivatives,
and $K(x,y)$ is the Gauss curvature, determined by $g_{ij}$ up
to their second derivatives.}

\medskip
This theorem holds even when $h_{ij}\in L^p$
({\it cf.} Maradare \cite{Mardare1,Mardare2}).
Note that system \eqref{g1} with \eqref{g2} is a system of nonlinear PDEs
of mixed elliptic-hyperbolic type, which is determined by the sign of the Gauss curvature
$K$.
Surfaces with Gauss curvature of changing sign
are very normal in geometry, including
tori such as
toroidal shells or doughnut surfaces.

\medskip
{\bf Fluid dynamics formalism for isometric embedding} (Chen-Slemrod-Wang \cite{CSW10a}).
Set $L=\r v^2+p, M=-\r uv, N=\r u^2+p$, and $q^2=u^2+v^2$.
Choose  $p$ as  the Chaplygin type gas: $p=-\frac{1}{\r}.$

\smallskip
The Codazzi equations \eqref{g1} become the balance laws of
momentum equations:
$$
\left\{\begin{aligned} &\d_x(\r uv)+\d_y(\r v^2+p)
 =-(\r v^2+p)\G^{(2)}_{22}-2\r uv\G^{(2)}_{12}-(\r u^2+p)\G^{(2)}_{11}, \\[1mm]
&\d_x(\r u^2+p)+\d_y(\r uv)
 =-(\r v^2+p)\G^{(1)}_{22}-2\r uv\G^{(1)}_{12}-(\r u^2+p)\G^{(1)}_{11},
\end{aligned}\right.
$$
and the Gauss equation becomes the Bernoulli relation:
$$
p=-\sqrt{q^2+ K}.
$$
Define the sound speed: $c^2=p'(\r)$. Then $c^2=\frac{1}{\r^2}=q^2+K$.

\smallskip
$c^2>q^2$ and the ``flow" is subsonic  when
$K>0$;

\smallskip
$c^2<q^2$ and the ``flow" is supersonic  when $K<0$;

\smallskip
$c^2=q^2$ and the ``flow" is sonic when $K=0$.

\smallskip
Based on this connection, the existence and continuity
of isometric embeddings
via compensated compactness and entropy analysis
were first addressed in
Chen-Slemrod-Wang \cite{CSW10a}.

\medskip
For higher dimensional case,
the isometric embeddings of
$n$-dimensional Riemannian manifolds ($n\ge 3$) into $\R^d$
are described by the following Gauss-Codazzi-Ricci system:

\medskip
\smallskip
\noindent
{\bf Gauss equations}:
\begin{equation}\label{GE}
h_{ji}^a h_{kl}^a-h_{ki}^a h_{jl}^a=R_{ijkl};
\end{equation}
{\bf Codazzi equations}:
\begin{equation}\label{CE}
\frac{\d h_{lj}^a}{\d x^k}-\frac{\d h_{kj}^a}{\d x^l} +
\Gamma_{lj}^mh_{km}^a-\Gamma_{kj}^m h_{lm}^a +\kappa_{kb}^a
h_{lj}^b-\kappa_{lb}^ah_{kj}^b=0;
\end{equation}
{\bf Ricci equations}:
\begin{equation}\label{RE}
\frac{\d\kappa_{lb}^a}{\d x^k}-\frac{\d\kappa_{kb}^a}{\d x^l}
-g^{mn}\left(h^a_{ml}h^b_{kn}-h^a_{mk}h^b_{ln}\right)
+\kappa_{kc}^a\kappa_{lb}^c-\kappa_{lc}^a\kappa_{kb}^c=0,
\end{equation}
where
$\{R_{ijkl}\}$ is the Riemann curvature tensor,
$\kappa_{kb}^a=-\kappa_{ka}^b$ are the coefficients of the
connection form (torsion coefficients) on the normal bundle;
the indices $a, b, c$ run from $1$ to $N$, and $i, j, k, l, m, n$ run
from $1$ to $d\ge 3$.

The Gauss-Codazzi-Ricci system \eqref{GE}--\eqref{RE} has no
type, neither purely hyperbolic nor purely elliptic for general
Riemann curvature tensor $R_{ijkl}$;
see Bryant-Griffiths-Yang \cite{BGY}.
Even though, we have established the following weak continuity and rigidity
of system \eqref{GE}--\eqref{RE} and corresponding embedded surfaces:

\medskip

{\bf Theorem} (Chen-Slemrod-Wang \cite{CSW10b}).  {\it Consider
the Gauss-Codazzi-Ricci system \eqref{GE}--\eqref{RE}.
\begin{enumerate}\renewcommand{\theenumi}{\roman{enumi}}
\item[\rm (i)] Let $(h_{ij}^{a,\varepsilon}, \kappa_{lb}^{a,\varepsilon})$ be a
sequence of solutions to system \eqref{GE}--\eqref{RE}, which is
uniformly bounded in $L^p, p>2$. Then the weak limit vector
field $(h_{ij}^{a}, \kappa_{lb}^{a})$ of the sequence
$(h_{ij}^{a,\varepsilon}, \kappa_{lb}^{a,\varepsilon})$ in
$L^p$ is still a solution to system \eqref{GE}--\eqref{RE}.

\item[\rm (ii)] There exists a minimizer $(h_{ij}^a, \kappa_{lb}^a)$
for the minimization problem:
$$
\min_S \|(h,\kappa)\|_{L^p(\Omega)}^p:= \min_S
\int_{\Omega}\sqrt{|g|}\left(|h_{ij}h_{ij}|^{\frac{p}{2}}+
|\kappa_{lb}\kappa_{lb}|^{\frac{p}{2}}\right) dx,
$$
where $S$ is the set of weak solutions to
system \eqref{GE}--\eqref{RE}.
\end{enumerate}
}

This weak continuity and rigidity theorem
is a reminiscence of the polyconvexity theory
in nonlinear elasticity by Ball \cite{Ball}, for which the rigidity
of elastic bodies can be achieved.

The proof of this theorem is based on
the following observations on the div-curl structure of the
Gauss-Codazzi-Ricci system:
\begin{eqnarray*}
&& {\text{div}}\,(\underbrace{\overbrace{0, \cdots, 0,
h_{lj}^{a,\e}}^{k}, 0, \cdots, -h_{kj}^{a,\e}}_{l},0,\cdots,0)=R_1,\\
&&{\text{curl}}\,(h^{a,\e}_{1i}, h^{a,\e}_{2i}, \cdots,
h^{a,\e}_{di})=R_2,\\
&&{\text{div}}\,(\underbrace{\overbrace{0, \cdots, 0,
\kappa_{lb}^{a,\e}}^{k}, 0, \cdots, -\kappa_{kb}^{a,\e}}_l, 0,
\cdots, 0)=R_3, \\
&&{\text{curl}}\,(\k^{a,\e}_{1b},
\k^{a,\e}_{2b}, \cdots, \k^{a,\e}_{db})=R_4,\\
&&{\text{div}}\,(\underbrace{\overbrace{0, \cdots, 0,
h_{li}^{b,\e}}^{k}, 0, \cdots, -h_{ki}^{b,\e}}_{l},0,\cdots, 0)=R_5,\\
&&{\text{curl}}\,(h^{b,\e}_{1i}, h^{b,\e}_{2i}, \cdots,
h^{b,\e}_{di})=R_6,\\
&&{\text{div}}\,(\underbrace{\overbrace{0, \cdots, 0,
\kappa_{lc}^{b,\e}}^{k}, 0,\cdots, -\kappa_{kc}^{b,\e}}_l, 0,\cdots,
0)=R_7,\\
&&{\text{curl}}\,(\k^{b,\e}_{1c},\k^{b,\e}_{2c},\cdots,\k^{b,\e}_{dc})=R_8,
\end{eqnarray*}
and $R_j, j=1,2, \cdots, 8$,
are confined in a compact set in
$\, H^{-1}_{loc}(\Omega)$.

Then employing the Murat-Tartar's div-curl lemma directly yields
\begin{eqnarray*}
h_{lj}^{a,\e}h_{ki}^{a,\e}-h_{kj}^{a,\e}h_{li}^{a,\e}&\rightharpoonup&
h_{lj}^{a}h_{ki}^{a}-h_{kj}^{a}h_{li}^{a}, \label{3.9a}\\
h_{lj}^{a,\e}h_{ki}^{b,\e}-h_{kj}^{a,\e}h_{li}^{b,\e}&\rightharpoonup&
h_{lj}^{a}h_{ki}^{b}-h_{kj}^{a}h_{li}^{b},\label{3.9}\\
\kappa_{kb}^{a,\e}\kappa_{lc}^{b,\e}-\kappa_{lb}^{a,\e}\kappa_{kc}^{b,\e}
 &\rightharpoonup&
\kappa_{kb}^a\kappa_{lc}^b-\kappa_{lb}^a\kappa_{kc}^b,\label{3.10}\\
\kappa_{kb}^{a,\e}h_{li}^{b,\e}-\kappa_{lb}^{a,\e} h_{ki}^{b,\e}
&\rightharpoonup& \kappa_{kb}^ah_{li}^b-\kappa_{lb}^a
h_{ki}^b\label{3.11}
\end{eqnarray*}
in the sense of distributions as $\varepsilon\to 0$, which implies
the weak continuity and rigidity of system \eqref{GE}--\eqref{RE}
and corresponding isometric embeddings.

A compactness framework for the
Gauss-Codazzi-Ricci system \eqref{GE}--\eqref{RE} has also established
in \cite{CSW10b}: {\it Given any sequence of
approximate solutions to this system which is uniformly bounded in
$L^2$ and has reasonable bounds on the errors made in the
approximation (the errors are confined in a compact subset of
$H^{-1}_{\text{loc}}$), then the approximating sequence has a weakly
convergent subsequence whose limit is still a solution
of system \eqref{GE}--\eqref{RE}.}

These results indicate that
the weak limit of isometrically embedded surfaces is still an
isometrically embedded surface in $\R^d$ for any Riemann curvature
tensor $\{R_{ijkl}\}$ without restriction, which is the rigidity property
of embedded surfaces in geometry.

\bigskip
\noindent
{\bf Acknowledgements.}  The materials presented above
include direct and/or indirect contributions
of my collaborators Xiaxi Ding, Feimin Huang, Philippe LeFloch,
Bang-He Li, Tianhong Li, Yunguang Lu,
Peizhu Luo, Mikhail Perepelitsa,
Marshall Slemrod, Dehua Wang, Tian-Yi Wang, Yongqian Zhang, Dianwen Zhu,
among others.

% ----------------------------------------------------------------


\begin{thebibliography}{99}


\bibitem{Alberti-Muller}
G. Alberti and S. M\"{u}ller,
A new approach to variational problems with multiple scales.
Comm. Pure Appl. Math. 54 (2001), 761--825.

\bibitem{Ball}
J.~M. Ball,
Convexity conditions and existence theorems in nonlinear elasticity.
Arch. Rational Mech. Anal. 63 (1976/77), 337--403.

\bibitem{Ball2}
J.~M. Ball, {A version of the fundamental theorem for Young
measures}.
{Lecture Notes in Phys.} {344}, pp. 207--215, Springer:
Berlin, 1989.


\bibitem{BB} S. Bianchini and A. Bressan,  Vanishing viscosity solutions of
nonlinear hyperbolic systems.
Ann. of Math. (2), 161 (2005), 223--342.

\bibitem{Bressan-book}
A. Bressan, Hyperbolic Systems of Conservation Laws: The One-Dimensional Cauchy Problem.
Oxford University Press: Oxford, 2000.

\bibitem{BCM}
M. Briane, J. Casado-D¨ªaz, and F. Murat,
The div-curl lemma ``trente ans apr¨¨s":
an extension and an application to the G-convergence of unbounded monotone operators.
J. Math. Pures Appl. (9) 91 (2009), 476--494.


\bibitem{BGY}
R.  L. Bryant,  P. A. Griffiths, and D. Yang,
Characteristics and existence of isometric embeddings.
Duke Math. J. 50 (1983), 893--994.

\bibitem{Chen3} G.-Q. Chen,  Convergence of the Lax-Friedrichs
scheme for isentropic gas dynamics (III).
Acta Math. Sci. 6B (1986), 75--120 (in English);
8A (1988), 243--276 (in Chinese).


\bibitem{Chen97}  G.-Q. Chen,
Remarks on spherically symmetric solutions of the compressible Euler equations.
Proc. Royal Soc. Edinburgh, 127A (1997), 243--259.


\bibitem{Chen05} G.-Q. Chen,
 Euler Equations and Related Hyperbolic Conservation Laws.
 In: {Handbook of Differential Equations: Evolutionary
       Differential Equations}, Vol. {\bf 2}, pp. 1--104, 2005,
Eds. C.~M. Dafermos and E. Feireisl, Elsevier Science B.V:
Amsterdam, The Netherlands.

\bibitem{ChenChristoforou}
G.-Q. Chen and C. Christoforou,
Solutions for a nonlocal conservation law with fading memory.
Proc. Amer. Math. Soc. 135 (2007), no. 12, 3905--3915.

\bibitem{CDSW}
G.-Q. Chen, C.~M. Dafermos, M. Slemrod, M., and D. Wang, On
two-dimensional sonic-subsonic flow. {Commun. Math. Phys.}
{271} (2007), 635--647.

\bibitem{CHW14}
G.-Q. Chen, F.-M. Huang, and T.-Y. Wang,
Subsonic-sonic limit of approximate solutions to multidimensional
steady Euler equations.
{Arch, Rational Mech. Anal.} 2015 (to appear);
arXiv:1311.3985.

\bibitem{CLF} G.-Q. Chen and P.~G. LeFloch,
{Compressible Euler equations with general pressure law}, {Arch.
Rational Mech. Anal.} {153} (2000), 221--259; {Existence theory
for the isentropic Euler equations}, {Arch. Rational Mech.
Anal.} {\bf 166} (2003), 81--98.


\bibitem{ChenLi}  G.-Q. Chen and T.-H. Li,
Well-posedness for two-diemnsonal steady supersonic Euler flows
past a Lipschitz wedge.
J. Diff. Eqs. 244 (2008), 1521--1550.


\bibitem{ChenLiLi}  G.-Q. Chen, B.-H. Li, and T.-H. Li,
Entropy solutions in $L^\infty$ for the Euler equations in nonlinear elastodynamics and
related equations. Arch. Rational Mech. Anal. 170 (2003), 331--357.



\bibitem{ChenLu}  G.-Q. Chen and Y.-G. Lu,
The study on application way of the compensated compactness theory.
Chinese Sci. Bull. 33 (1988), 641--644 (in Chinese);
34 (1989), 15--19 (in English).


\bibitem{CP} G.-Q. Chen and M. Perepelitsa,
{ Vanishing viscosity limit of the Navier-Stokes equations to the Euler equations
for compressible fluid flow.\/}
Comm. Pure Appl. Math. 63 (2010), 1469--1504.

\bibitem{CP2} G.-Q. Chen and M. Perepelitsa,
Shallow water equations: viscous solutions and inviscid limit.
Z. Angew. Math. Phys. 63 (2012), 1067--1084.

\bibitem{CP3} G.-Q. Chen and M. Perepelitsa,
Vanishing viscosity solutions of the compressible Euler equations with
 spherically symmetry and large initial data.
{Comm. Math. Phys.} {338} (2015), 771--800.


\bibitem{CSW08}
G.-Q. Chen, M. Slemrod, and D. Wang, Vanishing viscosity method
for transonic flow. { Arch. Rational Mech. Anal.}  189 (2008), 159--188.

\bibitem{CSW10a}
G.-Q. Chen,  M. Slemrod, and D. Wang,
{ Isometric immersions and compensated compactness}.
Comm. Math. Phys. 294 (2010), 411--437.

\bibitem{CSW10b}
G.-Q. Chen, M. Slemrod, and D. Wang, Weak continuity of the
Gauss-Codazzi-Ricci system for isometric embedding.
Proc. Amer. Math. Soc. 138 (2010), 1843--1852.


\bibitem{CZZ}
G.-Q. Chen, Y.-Q. Zhang, and D.-W. Zhu,
Existence and stability of supersonic Euler flows past Lipschitz wedges.
Arch. Rational Mech. Anal. 181 (2006), 261--310.


\bibitem{CFr} R. Courant and  K.~O. Friedrichs,
{Supersonic Flow and Shock Waves}.
Springer-Verlag: New York, 1948.

\bibitem{Dafermos88} C. M. Dafermos, Solutions in $L^\infty$ for a conservation law with memory.
Analyse Math\'{e}matique et Applications, 117--128, Gauthier-Villars,
Montrouge, 1988.

\bibitem{Dafermos} C. M. Dafermos,
Hyperbolic Conservation Laws in Continuum Physics.
Springer-Verlag: Berlin, 2010.

\bibitem{DeLellis-Szekelyhidi-1}
C. De Lellis and L. Sz\'{e}kelyhidi, Jr.,
On admissibility criteria for weak solutions of the Euler equations.
Arch. Rational Mech. Anal. 195 (2010),  225--260.

\bibitem{DeLellis-Szekelyhidi-2}
C. De Lellis and L. Sz\'{e}kelyhidi, Jr.,
The h-principle and the equations of fluid dynamics.
Bull. Amer. Math. Soc. (N.S.) 49 (2012), 347--375.


\bibitem{DCL} X. Ding, G.-Q. Chen, and P. Luo,
Convergence of
the Lax-Friedrichs scheme for the isentropic gas dynamics (I)-(II),
Acta Math. Sci. {5B} (1985), 483--500, 501--540 (in
English); {7A} (1987), 467-480; {8A} (1989), 61--94 (in
Chinese);  Convergence of the fractional step Lax-Friedrichs
scheme and Godunov scheme for the isentropic system of gas
dynamics, {Comm. Math. Phys.} {121} (1989), 63--84.

\bibitem{DiPerna77}
R.~J. DiPerna, Decay of solutions of hyperbolic systems of
conservation laws with a convex extension.
Arch. Rational Mech. Anal. 64 (1977), 1--46.

\bibitem{DiPerna83a} R.~J. DiPerna,
{Convergence of the viscosity method for isentropic gas dynamics.\/}
{Commun. Math. Phys.} {91} (1983), 1--30.

\bibitem{DiPerna83b} R.~J. DiPerna,
Convergence of approximate solutions to conservation laws.
Arch. Rational
Mech. Anal. 82 (1983), 27--70.


\bibitem{DiPerna85a} R.~J. DiPerna,
Measure-valued solutions to conservation laws.
Arch. Rational Mech. Anal. 88 (1985), 223--270.

\bibitem{DiPerna85b}
R.~J. DiPerna, Compensated compactness and general systems of
conservation laws. {Trans. Amer. Math. Soc.} {292} (1985),
383--420.

\bibitem{Efimov-a}  N.~V. Efimov, The impossibility in Euclideam
3-space of a complete regular surface with a negative upper bound of
the Gaussian curvature. {Dokl. Akad. Nauk SSSR (N.S.)}, {\bf
150} (1963), 1206--1209; {Soviet Math. Dokl.} {\bf 4} (1963),
843--846.


\bibitem{evans} L.~C. Evans, {Weak Convergence Methods for Nonlinear
Partial Differential Equations}. CBMS-RCSM, {74}, AMS:
Providence, RI, 1990.

\bibitem{Gilbarg}
D. Gilbarg, The existence and limit behavior of the one-dimensional shock layer.
Amer. J. Math.
73 (1951), 256--274.



\bibitem{Glimm}
J. Glimm, Solutions in the large for nonlinear hyperbolic systems
of equations. Comm. Pure Appl. Math. 18 (1965), 697--715.

\bibitem{Glimm-Lax}
J. Glimm and P. D. Lax,
Decay of solutions of systems of hyperbolic conservation laws.
Bull. Amer. Math. Soc. 73 (1967), no. 105.


\bibitem{Gromov86}
M. Gromov, {Partial Differential Relations}. Springer-Verlag:
Berlin, 1986.

\bibitem{Guderley}
G. Guderley, Starke kugelige und zylindrische Verdichtungsstosse
inder Nahe des Kugelmittelpunktes bzw. der Zylinderachse.
Luftfahrtforschung 19 (1942), no. 9, 302--311.

\bibitem{GMWZ}
C.M.I.O. Gu\`{e}s, G. M\'{e}tivier, M. Williams, and K. Zumbrun,
Navier-Stokes regularization of
multidimensional Euler shocks. Ann. Sci. \'{E}cole Norm. Sup. (4) 39 (2006),
75--175.

\bibitem{Gunther89}
M. G\"{u}nther,
Zum Einbettungssatz von J. Nash. (German) [On the embedding theorem of J. Nash].
Math. Nachr. 144 (1989), 165--187.

\bibitem{HanHong}
Q. Han  and J.-X. Hong,  {Isometric Embedding of Riemannian Manifolds in Euclidean Spaces}.
AMS: Providence, RI, 2006.

\bibitem{Hoff}
D. Hoff, Global solutions of the equations of one-dimensional, compressible flow with large
data and forces, and with differing end states.
Z. Angew. Math. Phys. 49 (1998), 774--785.

\bibitem{Hoff-Liu}
D. Hoff and T.-P. Liu, The inviscid limit for the Navier-Stokes equations of compressible, isentropic
flow with shock data. Indiana Univ. Math. J. 38 (1989), 861--915.

\bibitem{Holden-Risebro}
H. Holden and N.~H. Risebro,
Front Tracking for Hyperbolic Conservation Laws.
Springer: New York, 2011.


\bibitem{Hopf}
E. Hopf, The partial differential equation $u_t + uu_x = ¦Ìu_{xx}$.
Comm. Pure Appl. Math. 3 (1950), 201--230.

\bibitem{Huang-Wang}
F. Huang and Z. Wang,
Convergence of viscosity solutions for isothermal gas dynamics.
SIAM J. Math. Anal. 34 (2002), 595--610.

\bibitem{Janet}
M. Janet, {Sur la possibilit\'e de plonger un espace riemannian donn\'e dans un espace euclidien}.
Ann. Soc. Pol. Math. 5 (1926), 38--43.



\bibitem{Kanel}
I. Kanel, On a model system of equations for one-dimensional gas motion.
Differ. Uravn. 4 (1968), 721--734.

\bibitem{Kruzkov}
S. N. Kruzkov,
First order quasilinear equations with several independent variables
(Russian). Mat. Sb. (N.S.) 81 (123) (1970), 228--255.


\bibitem{Laxone} P. D. Lax, \, Shock wave and entropy.
In: Contributions to Functional Analysis,
 ed. E.A. Zarantonello, 603--634, Academic Press: New York, 1971.

\bibitem{Lax73}
P. D. Lax, Hyperbolic Systems of Conservation Laws and the Mathematical Theory of Shock
Waves. CBMS Regional Conference Series in Mathematics, No. 11. Philadelphia:
SIAM, 1973.




\bibitem{LeFloch-book}
P.~G. LeFloch,
Hyperbolic Systems of Conservation Laws:
The Theory of Classical and Nonclassical Shock Waves.
Birkh\"{a}user Verlag: Basel, 2002.

\bibitem{LeFloch-Shelukhin}
P.~G. LeFloch and V. Shelukhin,
Symmetries and local solvability of the isothermal gas dynamics equations.
Arch. Rational Mech. Anal. 175 (2005), 389--430.


\bibitem{LW}
P.~G. LeFloch and M. Westdickenberg,  \, Finite energy solutions to
the isentropic Euler equations with geometric effects. {J.
Math. Pures Appl.} {88} (2007), 386--429.


\bibitem{Liu81}
T.-P. Liu, Admissible solutions of hyperbolic conservation laws.
Mem. Amer. Math. Soc. 30 (1981), no. 240.

\bibitem{LiuYang}
T.-P. Liu and T. Yang,
Well-posedness theory for hyperbolic conservation laws.
Comm. Pure Appl. Math. 52 (1999), 1553--1586.


\bibitem{LPS}  P.-L. Lions, B. Perthame, and P.~E. Souganidis,  \,
{Existence and stability of entropy solutions for the hyperbolic
systems of isentropic gas dynamics in Eulerian and Lagrangian
coordinates}. Comm. Pure Appl. Math. 49 (1996),
599--638.

\bibitem{LPT} P.-L. Lions, B. Perthame, and E. Tadmor, \,
 Kinetic formulation of the isentorpic gas
dynamics and p-systems. {Comm. Math. Phys.} {163} (1994), 415--431.


\bibitem{Mardare1} S. Mardare,
The fundamental theorem of surface theory for surfaces with little
regularity. {J. Elasticity}, {73} (2003), 251--290.

\bibitem{Mardare2} S. Mardare,
On Pfaff systems with $L^p$ coefficients and their applications in
differential geometry. {J. Math. Pure Appl.} {84} (2005),
1659--1692.


\bibitem{Mor85}
C.~S. Morawetz,  {On a weak solution for a transonic flow problem}.
{Comm. Pure Appl. Math.} {38} (1985), 797--818.

\bibitem{Morawetz91}
C.~S. Morawetz, An alternative proof of DiPerna's theorem.
Comm. Pure Appl. Math. 44 (1991), 1081--1090.

\bibitem{Mor95}
C.~S. Morawetz,
 {On steady transonic flow by compensated compactness.}
 {Methods Appl. Anal.} {2} (1995), 257--268.


\bibitem{Murat} F. Murat,  Compacit\'e par compensation.
{Ann. Scuola Norm. Sup. Pisa Sci. Fis. Mat.} {5} (1978),
489--507.

\bibitem{Murat79} F. Murat,
Compacit\'{e} par compensation. II (French).
In: Proceedings of the International Meeting on Recent Methods in Nonlinear Analysis (Rome, 1978),
pp. 245--256, Pitagora, Bologna, 1979.

\bibitem{Murat81a} F. Murat,
L'injection du c\^{o}ne positif de $H^{-1}$ dans $W^{-1,q}$ est compacte pour
tout $q<2$ (French).
J. Math. Pures Appl. (9) 60 (1981), 309--322.

\bibitem{Murat81b} F. Murat,
Compacit\'{e} par compensation: condition n\'{e}cessaire et suffisante de continuit\'{e} faible
sous une hypoth¨¨se de rang constant (French).
Ann. Scuola Norm. Sup. Pisa Cl. Sci. (4) 8 (1981), no. 1,
69--102.

\bibitem{Murat87} F. Murat,
A survey on compensated compactness. In: Contributions to Modern Calculus of Variations (Bologna, 1985),
145--183, Pitman Res. Notes Math. Ser. 148, Longman Sci. Tech., Harlow, 1987.


\bibitem{Nash56}
J. Nash,
\newblock The imbedding problem for Riemannian manifolds.
{Ann. Math. (2)}, {63} (1956), 20--63.

\bibitem{Oleinik}
O.~A. Oleinik,
Discontinuous solutions of non-linear differential equations. Usp. Mat. Nauk 12
(1957), 3--73. English translation: AMS Translations, Ser. II, 26, 95--172.

\bibitem{Perthame-Tzavaras}
B. Perthame and A.~E. Tzavaras,
Kinetic formulation for systems of two conservation laws and elastodynamics.
Arch. Rational Mech. Anal. 155 (2000), 1--48.

\bibitem{Rauch}
J. Rauch,
BV estimates fail for most quasilinear hyperbolic systems in dimension greater
than one. Comm. Math. Phys. 106 (1986), 481--484.

%R
\bibitem{Ro} S. Rosseland,
 {The Pulsation Theory of Variable Stars}.
    Dover Publications, Inc.: New York, 1964.

\bibitem{Schonbek}
M.~E. Schonbek,
Convergence of solutions to nonlinear dispersive equations.
Comm. Partial Diff. Eqs. 7 (1982), 959--1000.

\bibitem{Serre87}
D. Serre, La compacit\'{e} par compensation pour les syst\`{e}mes non lin\'{e}aires de deux equations
a une dimension d'espace.
J. Math. Pures Appl. 65 (1987), 423--468.


\bibitem{SS}
D. Serre and J.~W. Shearer,  Convergence with physical viscosity
for nonlinear elasticity. Preprint, 1994 (unpublished).

\bibitem{Sm}
J. Smoller,
Shock Waves and Reaction-Diffusion Equations.
Second edition. Springer-Verlag: New York, 1994.


\bibitem{TRB}
E. Tadmor, M. Rascle, and P. Bagnerini,
Compensated compactness for 2D conservation laws.
J. Hyper. Diff. Eqs. 2 (2005), 697--712.


\bibitem{Tartar}
L. Tartar,   Compensated compactness and
applications to partial differential equations. In: {Research
Notes in Mathematics, Nonlinear Analysis and Mechanics},
Herriot-Watt Symposium, Vol. {4}, Knops R.J. ed., Pitman Press,
1979.

\bibitem{Tartar83a}
L. Tartar,
The compensated compactness method applied to systems of conservation laws.
In: Systems of Nonlinear Partial Differential Equations (Oxford, 1982), 263--285,
NATO Adv. Sci. Inst. Ser. C Math. Phys. Sci. 111, Reidel, Dordrecht, 1983.

\bibitem{Tartar83b}
L. Tartar, Compacit\'{e} par compensation: r\'{e}sultats et perspectives (French).
 In: Nonlinear Partial Differential Equations and Their Applications,
 Coll\`{e}ge de France Seminar, Vol. IV (Paris, 1981/1982), 350--369,
 Res. Notes in Math. 84, Pitman, Boston, MA, 1983.

\bibitem{Tartar84}
L. Tartar,
Oscillations in nonlinear partial differential equations: compensated compactness and homogenization.
In: Nonlinear Systems of Partial Differential Equations in Applied Mathematics, Part 1 (Santa Fe, N.M., 1984),
243--266, Lectures in Appl. Math., 23, Amer. Math. Soc.: Providence, RI, 1986.

\bibitem{Tartar87}
L. Tartar, Discontinuities and oscillations.
In: Directions in Partial Differential Equations (Madison, WI, 1985), 211--233,
Publ. Math. Res. Center Univ. Wisconsin, 54, Academic Press, Boston, MA, 1987.

\bibitem{Tartar08}
L. Tartar, From Hyperbolic Systems to Kinetic Theory: A Personalized Quest.
Lecture Notes of the Unione Matematica Italiana, 6. Springer-Verlag, Berlin; UMI, Bologna, 2008.

\bibitem{Tartar09}
L. Tartar, The General Theory of Homogenization: A Personalized Introduction.
Lecture Notes of the Unione Matematica Italiana, 7.
Springer-Verlag, Berlin; UMI, Bologna, 2009.


\bibitem{Volpert}
A. I. Vol'pert, Spaces BV and quasilinear equations (Russian).
Mat. Sb. (N.S.) 73 (115) (1967), 255--302.


%Wh
\bibitem{Wh} G.~B. Whitham,
{Linear and Nonlinear Waves}.
    John Wiley \& Sons: New York, 1974.

\bibitem{Yau00}
S.-T. Yau, Review of geometry and analysis. In: {Mathematics:
Frontiers and Perspectives}, pp. 353--401, International Mathematics
Union, Eds. V. Arnold, M. Atiyah, P. Lax, and B. Mazur, AMS:
Providence, 2000.
\end{thebibliography}
\end{document}